\documentclass[12pt,reqno]{amsart}

\usepackage{amssymb, amsmath,}

\newtheorem{theorem}{Theorem}[section]

\newtheorem{corollary}{Corollary}[section]

\newtheorem{remark}{Remark}[section]

\begin{document}

\title{  Bohr and Rogosinski abscissas  for ordinary Dirichlet  series. }
\author {
Lev Aizenberg, Victor Gotlib, Alekos Vidras}
\thanks{AMS classification number: 30B50, 30A10.}
\thanks{Keywords: Isometric Bohr abscissa, Rogosinski abscissa,
 Dirichlet series.}
 \address{Department of Mathematics, Bar-Ilan University\newline
 Ramat-Gan 52900, Israel.}
\address{Department of Mathematics, Holon Institute of Technology,\newline
 P.O.B  305, Holon 58102, Israel.}
\address{Department of Mathematics and Statistics,\newline Univ.of Cyprus,
Nicosia 1678, Cyprus.}
 \email{aizenbrg@macs.biu.ac.il, gotlib@hit.ac.il}
 \email{msvidras@ucy.ac.cy}
 \maketitle
 \begin{abstract}
 We prove that the  abscissas  of Bohr and Rogosinski for ordinary Dirichlet
 series, mapping the right half-plane into the bounded convex
 domain $G\subset \mathbb{C} $ are independent of the domain $G$.
 Furthermore, we obtain new estimates about these  abscissas.
 \end{abstract}
 \section{Preliminaries }
 \setcounter{equation}{0}
 Let us recall the theorem of H.Bohr \cite{bohr3:gnus} in 1914.
\begin{theorem} If a power series
\begin{eqnarray}
f(z_1)=\sum\limits_{k=0}^{\infty}c_kz_1^k
\end{eqnarray}
 converges in the
unit disk $U_1$ and its sum has modulus
 less than $1$, then
 \begin{eqnarray}
 \sum\limits_{k=0}^{\infty }\vert c_k z_1^k\vert <1,
 \end{eqnarray}
 if $ \vert z_1\vert <{1\over 3}$. Moreover,
 the constant ${1\over 3}$ cannot be improved.
 \end{theorem}
 For reasons of convenience we write the inequality $(1.2)$ in the
 following equivalent form
 \begin{eqnarray*}
\sum\limits_{k=1}^{\infty }\vert c_k z_1^k\vert <1-\vert c_0\vert .
\end{eqnarray*}
Generalizations and modifications  of this result can be found in
\cite{aiz2:gnus}, \cite{aiz4:gnus}, \cite{aizad1:gnus} --
\cite{aizg:gnus}, \cite{aizvi:gnus}, \cite{ben:gnus},
{\cite{dra:gnus}, \cite{kres:gnus},  \cite{pau:gnus},
\cite{sid:gnus}, \cite{tom:gnus} . On the other hand, formulation
of Bohr theorem in several complex variables appeared very
recently. Given a complete Reinhardt domain ${\mathcal{D}}\subset
\mathbb{C}^n$, we denote by $R( {\mathcal D})$ the largest
non-negative number $r$ with the property that if the power series
 \begin{eqnarray}
 f(z)=\sum\limits_{\vert \alpha \vert \geq 0}c_{\alpha}z^{\alpha},
 \; z\in {\mathcal D},
 \end{eqnarray}
 where as usual, $\alpha=(\alpha_1,\dots, \alpha_n)\in \mathbb{N}^n_0$,
 $z^{\alpha}=z_1^{\alpha_1}\dots z_n^{\alpha_n}$, $\vert \alpha
 \vert =\alpha_1+\dots \alpha_n $, converges in ${\mathcal D}$ and the
 modulus of its sum is less than $1$, then
 \begin{eqnarray*}
 \sum\limits_{\vert \alpha \vert \geq 1}\vert c_{\alpha}z^{\alpha}
 \vert < 1-\vert c_0\vert
 \end{eqnarray*}
 in the homothetic domain $\mathcal{D}_r=r\mathcal{D}$. Here
 $c_0=c_{0,\dots, 0}$. A variety of results, related to this
 particular $R(\mathcal{D})$ or other multidimensional generalizations of Bohr radius are found in \cite{aiz1:gnus} - \cite{aizvi:gnus}, \cite{ben:gnus} - \cite{ boha:gnus},
  \cite{deff:gnus} -- \cite{dra:gnus}.
 Key results, used by us in the present paper, were obtained in
 \cite{aiz4:gnus}, \cite{aize:gnus}, \cite{bal:gnus}.\\
 Let $\widetilde G $ be a convex hull of the domain $G\subset \mathbb{C}$. For a a complete Reinhardt domain  $\mathcal{D}$ in $\mathbb{C}^n$, denote by
 $R(\mathcal{D},G)$
 the largest $r\geq 0 $ such that if the function $f(z)$ from $(1.3)$ is
 holomorphic in $\mathcal{D}$, $f(\mathcal{D})\subset G$ and  $\widetilde G \not=\mathbb{C}$,
 then
 \begin{eqnarray*}
 \sum\limits_{\vert \alpha \vert \geq 1}\vert
 c_{\alpha}z^{\alpha}\vert <dist (c_0,\partial\widetilde G)
 \end{eqnarray*}
in homothety $\mathcal{D}_r$. A point $p\in \partial G $ is called
{\it{a point of convexity}} if $p\in \partial\widetilde G$. A
point of convexity $p$ is called {\it {regular}} if there exists a
disk $U\subset G $
so that $p\in\partial U$. \\
The following result is a consequence of a more general result,
proved \cite{aiz4:gnus}. We state it in a suitable for us form for
the purposes of the present article.
\begin{theorem}
If $\widetilde G\not=\mathbb{C}$, then $R(\mathcal{D},G)$ is not
smaller than $R(\mathcal{D},U_1)$. If $\partial G $ contains at
least one point of regular convexity , then
\begin{eqnarray*}
R(\mathcal{D},G)=R(\mathcal{D},U_1).
\end{eqnarray*}
\end{theorem}
\begin{corollary}
If the domain $G$ is convex and $G\not= \mathbb{C}$, then
$R(\mathcal{D}, G)$ is independent of the choice of $G$.
\end{corollary}
Besides the Bohr radius we will use the radius of Rogosinski, whose
classic result of 1923, \cite{ rog:gnus, szego:gnus,laga:gnus}, is
described in the following statement:
\begin{theorem}
If the function $f(z_1)$ from $(1.1)$ is holomorphic in the unit
disk $U_1$ and $\vert f(z_1)\vert <1 $ in $U_1$, then all of its
partial sums are less than $1$ in the disk of radius ${1\over 2}$,
that is:
\begin{eqnarray*}
\vert \sum\limits_{k=0}^{m}c_kz_1^k\vert <1
\end{eqnarray*}
for $\vert z_1\vert <{1\over 2}$ and this radius is sharp.
\end{theorem}
The following result is the consequence of more general fact
proved in \cite{aize:gnus}. Again, we are going to use its
particular, convenient for us, formulation: let $\mathcal{A}$ be a
lattice in $\mathbb{N}_0^n$, which is represented by
\begin{eqnarray*}
\mathcal{A}=\{\alpha \in \mathbb{N}_0^n:\;
m_1\alpha_1+\dots+m_n\alpha_n\leq m\},
\end{eqnarray*}
where all the numbers $m_1,\dots, m_n, m$ belong to $\mathbb{N}_0$
and have no common divisor.
\begin{theorem}
Let the function $f(z)$ in $(1.3)$ be holomorphic in the Reinhardt
domain $\mathcal{D}$ and $f(\mathcal{D})\subset G$, where $G$ is a
convex domain in $\mathbb{C}$ so that $G\not=\mathbb{C}$. Then the
polynomial (partial sum)
\begin{eqnarray*}
\sum\limits_{\alpha \in \mathcal{ A}}c_{\alpha}z^{\alpha}
\end{eqnarray*}
maps $\mathcal{D}_{\mathcal{A}}$ into $G$, where
\begin{eqnarray*}
\mathcal{D}_{\mathcal{A}}=\{z\in\mathbb{C}^n:\;({z_1\over
r_m^{m_1}},\dots,{z_n\over r_m^{m_n}})\in \mathcal{D}\},
\end{eqnarray*}
does not depend on $G$, $r_1={1\over 2}$, $r_2={\sqrt 3\over 8}$
and for $l\geq 3$ the number  $r_l$ is the unique positive
solution of the equation
\begin{eqnarray*}
1-r-2r^{l+1}=0.
\end{eqnarray*}
\end{theorem}
In the present article we investigate the Dirichlet series
\begin{eqnarray}
f(s)=\sum\limits_{n=1}^{\infty}{a_n\over n^s},
\end{eqnarray}
converging in the right half-plane $\Pi=\{s\in \mathbb{C}:\;\Re s
=\sigma >0\}$ and $f(s)\in H^{\infty}(\Pi )$, that is
$\|f\|=\sup\limits_{s\in \Pi}\vert f(s)\vert <\infty $. Following
\cite{bal:gnus}, we call {\it{isometric Bohr abscissa}} ${\bf{b}}$
the non-negative real number defined as the infimum of those
$\sigma \geq 0 $ such that for all $f\in H^{\infty}(\Pi) $ which
can be expressed as a Dirichlet series $(1.4)$, the following
holds
\begin{eqnarray}
\sum\limits_{n=1}^{\infty}{\vert a_n\vert \over n^{\sigma}}<\|f\|.
\end{eqnarray}
 For reasons of convenience, we write $(1.5)$ as
 \begin{eqnarray*}
\sum\limits_{n=2}^{\infty} {\vert a_n\vert \over n^{\sigma}}\leq
\|f\|-a_1.
\end{eqnarray*}
The next result was obtained in \cite{bal:gnus}.
\begin{theorem}
For the isometric Bohr abscissa the following estimates are valid
\begin{eqnarray*}
1.5850\dots={\log 3\over \log 2}\leq {\bf{b}}<1.8154
\end{eqnarray*}
If $f\in H^{\infty}(\Pi )$ and $\|f\| =1 $, then
\begin{eqnarray*}
\vert a_1\vert ^2+\sum\limits_{n=2}^{\infty} {\vert a_n\vert \over
n^{\sigma}}\leq 1,
\end{eqnarray*}
where $\sigma <1.7287 $.
\end{theorem}
 Analogously, to the Bohr radius
$R(\mathcal{D},G)$, {\it{the isometric Bohr abscissa}}
${\bf{b}}(G)$, where $G$ is a domain in $\mathbb{C}$, $\widetilde
G\not=\mathbb{C}$, is defined as the infimum of those $\sigma \geq
0 $ such that if $f\in H^{\infty}(\Pi)$, $f(s)$ is like in $(1.4)$
and $f(\Pi)\subset G$, then
\begin{eqnarray}
\sum\limits_{n=2}^{\infty} {\vert a_n\vert \over
n^{\sigma}}<dist(a_1,\partial \widetilde G).
\end{eqnarray}
Furthermore, we define the {\it{ Rogosinski abscissa
${\bf{r}}(G)$}} for functions $f\in H^{\infty}(\Pi)$,
$f(\Pi)\subset G$, as the infimum of $\sigma ^{\prime}\geq 0 $,
such that for all partial sums of the series $(1.4)$
\begin{eqnarray}
\sum\limits_{n=1}^{k}{a_n\over n^s}=P_k(s)
\end{eqnarray}
the inclusion $P_k(\Pi _{\sigma ^{\prime}})\subset G$ is valid,
where $\Pi_{\sigma^{\prime}}=\{s\in \mathbb{C}: \;\sigma
>\sigma^{\prime}\}$.
\section{The main results}
\setcounter{equation}{0}
\begin{theorem}
Let $G$ be a be a bounded domain in $\mathbb{C}$. Then the
isometric Bohr abscissa ${\bf{b}}(G)$ is not larger then the
isometric Bohr abscissa ${\bf{b}}(U_1)$. If $\partial G$ contains
at least one regular point of convexity, then
${\bf{b}}(G)={\bf{b}}(U_1)$.
\end{theorem}
{\bf{Proof:}} We are going to exploit the connection between
classical Dirichlet series and the power series in the infinite
dimensional polydisc found by H.Bohr in \cite{bohr2:gnus}.
Consider the series $(1.4)$. Every natural number $n$ is product
of its prime factors
\begin{eqnarray*}
n=p_1^{\alpha_1}\dots p_m^{\alpha_m}.
\end{eqnarray*}
We set $z=(p_1^{-s}, p_2^{-s},\dots )$. Then
\begin{eqnarray}
f(s)=\sum\limits_{n=1}^{\infty}a_n(p_1^{-s})^{\alpha_1}\dots
(p_m^{-s})^{\alpha_m}=\sum\limits_{n=1}^{\infty}a_nz_1^{\alpha
_1}\dots z_m^{\alpha _m},
\end{eqnarray}
where $m=m(n)$. If the power series $(1.4)$ converges absolutely
for $\sigma =\sigma_0$, then the set of values of $f(s)$ on the
vertical line $\{s\in \mathbb{C}:\;\Re s=\sigma _0\}$ is
everywhere dense in the set of the values of the series $(2.1)$ on
the set $\{z:\; \vert z_1\vert ={1\over
p_1^{\sigma_0}},\dots,\vert z_m\vert ={1\over p_m^{\sigma_0}}\}$,
\cite{bohr2:gnus}. We will be using this result in the particular
case when the series $(1.4)$ is just a Dirichlet polynomial. The
same conclusion is easily deduced from a theorem of Kronecker,
\cite{bay:gnus}. Denote by $\sigma _u$ the abscissa of uniform
convergence of the series $(1.4)$, that is, the infimum of those
$\sigma $, so that the series $(1.4)$ converges uniformly in the
half-plane $\Pi_{\sigma}=\{s\in \mathbb{C}:\;\Re s>\sigma \}$. Let
$\sigma _{b}$ be infimum of those $\sigma $ for which the series
$(1.4)$ is bounded (and naturally holomorphic) in the half-plane
$\Pi_{\sigma}$.  A non-trivial result of Bohr states that
$\sigma_u=\sigma_{b}$, \cite{bohr1:gnus},\cite{bay:gnus}.
Therefore, if $f(\Pi)\subset G$, where $G$ is a bounded domain in
$\mathbb{C}$, then in every half-plane $\Pi_{\sigma}$, $\sigma >0$
the series $(1.4)$ converges uniformly. Hence, for $\sigma
$-fixed, we can find for every $\epsilon >0$ an index $k\in
\mathbb{ N}$ so that $\vert f(s)-P_k(s)\vert <\epsilon $, where
the polynomial $P_k(s)$ is taken from $(1.7)$.We remark that
$P_k(\Pi_{\sigma})\subset G_{\epsilon}$, where $G_{\epsilon}$ is
an $\epsilon $-neighborhood of $G$. Therefore the polynomial
\begin{eqnarray*}
\widetilde P_k(z)=\sum\limits_{n=1}^ka_nz_1^{\alpha_1}\dots
z_m^{\alpha_m}
\end{eqnarray*}
in the corresponding polydisc $\mathcal{D}$ is such that
$\widetilde P_k(\mathcal{D})\subset G_{\epsilon }$. From Theorem
1.2 then follows  that $R(\mathcal{D},G_{\epsilon})\geq
R(\mathcal{D}, U_1)$ and therefore the claim
\begin{eqnarray}
\sum\limits_{n=2}^k \vert a_n\vert {1\over p_1^{\sigma}}\dots
{1\over p_m^{\sigma}}<dist(a_1,\partial\widetilde G_{\epsilon})
\end{eqnarray}
for $\sigma \leq {\bf{b}}(U_1)$ is valid. Taking the limit in
$(2.2)$ for $k\longrightarrow \infty$ (this means $\epsilon
\longrightarrow 0$), we obtain that for such  a $\sigma $ the
relation $(1.6)$ holds, but instead of the proper inequality $< $,
the inequality $\leq $ is now possible. But if we had equality
there, then by taking smaller $\sigma $ we would have had obtained
the converse inequality, which is impossible. Therefore
${\bf{b}}(G)\leq {\bf{b}}(U_1)$.\par
 Let us assume now that on
the boundary $\partial G$ there exists at least one regular point
of convexity. Denote it by $p_0$. Then there exists a disk
$U\subset G$ so that $p_0\in\partial U\cap\partial G\cap \partial
\widetilde G$. Remark that homothety and  parallel translation for
the disk $U_1$ (that is all $a_n$ are replaced by $ra_n$ and after
this $ra_1$ is replaced by $ra_1+c$, and therefore the new disc
$U$ is $rU_1+ c$) do not alter the condition $(1.6)$ for the disc.
Hence ${\bf{b}}(U)={\bf{b}}(U_1)$. From here, it follows that
${\bf{b}}(G)\geq {\bf{b}}(U_1)={\bf{b}}(U)$ since the
corresponding set of holomorphic functions, which are represented
by Dirichlet series in $\Pi $ and are satisfying
$f(\Pi) \subset G$, is larger than the set of functions satisfying $f(\Pi)\subset U$. Thus, in this case we obtain ${\bf{b}}(G)={\bf{b}}(U_1)$.$\diamondsuit $.\\
\begin{corollary}
If $G$ is a bounded, convex domain in $\mathbb{C}$ then the
isometric Bohr abscissa ${\bf{b}}(G)$ does not depend from the
choice of the domain $G$.
\end{corollary}
We point out that in the proof of the Theorem 2.1, we used the
boundedness of the domain $G$ only to be able to approximate the
Dirichlet series by a Dirichlet polynomial, in order to use the
later's  uniform convergence in the corresponding half-plane and
thus to be able to apply Bohr result  about the everywhere density
of its values. Instead, one could demand that the series $(1.4)$
converges absolutely in the half plane $\Pi $. Thus, one is able
to formulate the following statement, where instead of isometric
Bohr abscissa ${\bf{b}}(G)$ one considers ${\bf{b}}_{a}(G)$,
defined in the same manner as ${\bf{b}}(G)$, but only for
absolutely converging Dirichlet series in the half plane $\Pi$.
\begin{theorem}
Let $G$ be a domain in $\mathbb{C}$, $\widetilde G
\not=\mathbb{C}$. Then the isometric Bohr abscissa
${\bf{b}}_{a}(G)$ is not greater ${\bf{b}}(U_1)$. If $\partial G$
contain at least one regular point of convexity, then
${\bf{b}}_{a}(G)={\bf{b}}(U_1)$.
\end{theorem}
\begin{theorem}
Let $G$ be a bounded convex domain in $\mathbb{C}$. Then the
following Rogosinski abscissas are equal: ${\bf {r}}(G)={\bf
{r}}(U_1)$.
\end{theorem}
{\bf{Proof:}} The proof  repeats the steps of the proof of the
Theorem 2.1, but instead of Theorem 1.2 one uses the Theorem 1.4
and the following remark: the Dirichlet polynomial $(1.7)$ can be
described via the inequality
\begin{eqnarray}
\alpha_1\log 2+\alpha_2\log 3+\dots +\alpha_m\log m\leq \log k,
\end{eqnarray}
where the prime numbers are all prime numbers that appear in the
prime decomposition of the integers $2,3,4,\dots, k$.
Approximating the logarithms in $(2.3)$ with fractions so that no
new integers appear and the existing ones remain, one can obtain
the lattice
\begin{eqnarray}
\alpha_1d+\alpha_2d_2+\dots +\alpha_md_m\leq d,
\end{eqnarray}
instead of the lattice $(2.3)$, where $d, d_i\in \mathbb{Q}$,
$i=1,\dots, m $. Crucial fact here is that the new lattice is
described by the same $m$-tuples $(\alpha_1,\dots, \alpha _m)$ as
in $(2.3)$. For the final step it is enough to obtain a common
denominator for the rational numbers $d, d_i$, $i=1,\dots, m $ and
get the lattice $\mathcal {A}$ from the Theorem
1.4.$\diamondsuit$\par
 We remark that defining the Rogosinski abscissa ${\bf {r}}_a(G)$ in
the same manner as ${\bf {r}}(G)$, but only for absolutely
convergent Dirichlet series in the half-plane, one obtains
\begin{theorem}
If the domain $G$ is convex then ${\bf {r}}_a(G)={\bf {r}}(U_1)$.
\end{theorem}
It is a well know fact that the Rogosinski radius is not smaller
than the Bohr radius for power series in one and several complex
variables (that is, the Rogosinski condition is satisfied in a
disc (ball) of larger radius than the Bohr condition). Therefore,
one might expect that the Rogosinski abscissa is not greater that
the isometric Bohr abscissa for ordinary Dirichlet series. This is
the content of the next result.
\begin{theorem}
Let $G$ be a convex bounded domain in $\mathbb{C}$. Then ${\bf
{r}}(G)\leq {\bf {b}}(G)$. If the domain $G\not=\mathbb{C}$ is
just convex, then ${\bf {r}}_a(G)\leq {\bf {b}}_a(G)$.
\end{theorem}
{\bf{Proof:}} Actually, if the domain $G$ is bounded  and $\sigma
>{\bf{b}}(G)$, then
\begin{eqnarray*}
\vert P_k(s)-a_1\vert =\vert \sum\limits_{n=2}^k{a_n\over
n^s}\vert \leq  \sum\limits_{n=2}^k\vert {a_n\over n^s}\vert\leq
\sum\limits_{n=2}^{\infty}\vert {a_n\over n^s}\vert
<dist(a_1,\partial G).
\end{eqnarray*}
But the obtained inequality
\begin{eqnarray*}
\vert P_k(s)-a_1\vert \leq dist(a_1,\partial G)
\end{eqnarray*}
means geometrically that $P_k(\Pi _{\sigma})\subset G$. Thus, if
$\sigma >{\bf {b}}(G)$, then $\sigma >{\bf {r}}(G)$ also.
Therefore ${\bf {r}}(G)\leq {\bf {b}}(G)$. The case when the
boundedness of the convex domain $G$ is not required is considered
analogously. $\diamondsuit $.
\begin{remark}
We point out the following open problems:\\
1. Is it possible to remove the condition on the boundedness of
the domain $G$ in the statements of the Theorems 2.1, 2.4 and the
Corollary 2.2?\\
2) Can one prove the second parts of the  Theorems 2.1, 2.3
without the assumption on the existence of at least one regular
point of convexity?
\end{remark}
\section{Estimates for the Bohr and Rogosinski abscissas}
\setcounter{equation}{0}
\begin{theorem}
Let $G$ be a convex bounded domain in $\mathbb{C}$. Then the
isometric Bohr abscissa satisfies
\begin{eqnarray}
{\bf{b}}(G)\leq 1.7267.
\end{eqnarray}
If $G$ is a domain in $\mathbb{C}$, $\widetilde G\not=
\mathbb{C}$, then the isometric Bohr abscissa ${\bf{b}}_a(G)$
satisfies the same estimate.
\end{theorem}
{\bf{Proof:}} Denote by $\Omega (n)$ the number of the prime
divisors (counted with their multiplicity ) of the natural number
$n$. Pivotal for us is the following result (\cite{bal:gnus},
Prop.2.1): let $f$ be like in $(1.4)$ and $\|f\|=1 $, then for
$k\geq 1$ one has
\begin{eqnarray}
(\sum\limits_{\Omega (n)=k, n\geq 2}\vert a_n\vert^2)^{{1\over
2}}\leq 1-\vert a_1\vert ^2.
\end{eqnarray}
Furthermore
\begin{eqnarray*}
\sum\limits_{k=1}^{\infty}{\vert a_k\vert\over n^{\sigma}}&\leq
&\vert a_1\vert+\sum\limits_{k=1}^{\infty}(\sum\limits_{{\Omega
(n)=k,\atop n\geq 2}}\vert a_n\vert^2)^{{1\over
2}}(\sum\limits_{{\Omega (n)=k,\atop n\geq 2}}{1\over
n^{2\sigma}})^{{1\over
2}}\\
&\leq & \vert a_1\vert+ (1-\vert
a_1\vert^2)\sum\limits_{k=1}^{\infty}(\sum\limits_{{\Omega
(n)=k,\atop n\geq 2}}{1\over n^{2\sigma}})^{{1\over
2}}\\
\end{eqnarray*}
Consider now the equation
\begin{eqnarray}
\sum\limits_{k=1}^{\infty}(\sum\limits_{{\Omega (n)=k,\atop n\geq
2}}{1\over n^{2\sigma}})^{{1\over 2}}={1\over 2}.
\end{eqnarray}
If $\sigma _0$ is the unique solution of the equation $(3.3)$,
then
\begin{eqnarray*}
\sum\limits_{k=1}^{\infty}{\vert a_k\vert\over n^{\sigma _0}}\leq
\vert a_1\vert+ (1-\vert a_1\vert ^2){1\over 2}\leq 1,
\end{eqnarray*}
since $\vert a_1\vert \leq 1$. Therefore ${\bf{b}}(U_1)\leq \sigma
_0$. The solution of the equation $(3.3)$ was obtained
numerically, using Maple. Thus the estimate $\sigma_0<1.7267$ was
obtained. $\diamondsuit $\\
\begin{corollary}
Let $G$ be a bounded domain in $\mathbb{C}$. Assume also that
$\partial G$ contains at least one regular point of convexity.
Then
\begin{eqnarray*}
1.5850\dots ={\log 3\over \log 2}\leq {\bf{b}}(G)<1.7267.
\end{eqnarray*}
If one does not require the boundedness of the domain
$G\not=\mathbb{C}$, then ${\bf{b}}_a(G)$ satisfies the same
estimates.
\end{corollary}
\begin{theorem}
Let $G$ be a convex bounded domain in $\mathbb{C}$. Then, if the
function $f$ is from $(1.5)$ and $\|f\|=1$, then for every such
$f$
\begin{eqnarray}
\vert a_1\vert^2+\sum\limits_{k=2}^{\infty}{\vert a_k\vert\over
n^{\sigma }}\leq 1,
\end{eqnarray}
where $1\leq \sigma < 1.2061$.
\end{theorem}
{\bf{Proof:}} From $(3.2)$, we obtain as before, that the left
hand-side in $(3.4)$ is not larger than
\begin{eqnarray*}
\vert a_1\vert ^2+ (1-\vert a_1\vert
^2)\sum\limits_{k=1}^{\infty}(\sum\limits_{{\Omega (n)=k,\atop
n\geq 2}}{1\over n^{2\sigma}})^{{1\over 2}}.
\end{eqnarray*}
Consider the equation
\begin{eqnarray}
\sum\limits_{k=1}^{\infty}(\sum\limits_{{\Omega (n)=k,\atop n\geq
2}}{1\over n^{2\sigma}})^{{1\over 2}}=1.
\end{eqnarray}
If $\sigma _0 $ is the root of the equation $(3.5)$, then for
$\sigma =\sigma _0$ the relation $(3.4)$ holds. The equation
$(3.5)$ was solved numerically by using Maple to obtain the
estimate
$\sigma_0<1.2061$.\\
On the other hand, a particular case of the Dirichlet series
$(1.4)$ is the series
\begin{eqnarray}
f(s)=\sum\limits_{n=0}^{\infty} {a_n\over 2^{ns}},
\end{eqnarray}
which is a power series relatively to the variable ${1\over 2^s}$.
It is known (\cite{pau:gnus}) that the best value for the  radius
$r$, for which the relation
\begin{eqnarray*}
\vert c_0\vert ^2+\sum\limits_{k=1}^{\infty}\vert c_k\vert r^k
\leq 1
\end{eqnarray*}
is valid for every power series $(1.1)$ satisfying $\vert
f(z_1)\vert <1 $ in the unit disk $U_1$, is ${1\over 2}$. Then for
the Dirichlet series $(3.6)$ we obtain
\begin{eqnarray*}
{1\over 2^{\sigma_0}}={1\over 2},
\end{eqnarray*}
that is $\sigma _0=1 $. Since the series $(3.6)$ is a particular
case of the series $(1.4)$, we deduce that $\sigma $ from Theorem
3.3 is greater or equal to 1.$\diamondsuit$\\
Analogously one can prove the following
\begin{theorem}
Let $G$ be a convex bounded domain in $\mathbb{C}$. Then the
Rogosinski abscissa ${\bf{r}}(G)\geq 1 $. If the domain $G\not
=\mathbb{C}$ is not bounded then the Rogosinski abscissa
${\bf{r}}_{a}(G)$ satisfies the same inequality.
\end{theorem}
\begin{remark}
We conclude the present paper by by stating the following
hypothesis: for every convex bounded domain $G\subset \mathbb{C}$
the equality ${\bf{r}}(G)=1 $ is true. If the domain $G\not
=\mathbb{C}$ is convex, but not necessarily bounded then also
${\bf{r}}_{a}(G)=1$.
\end{remark}

\end{document}